\documentclass[11pt,leqno]{amsart}
\usepackage{amsmath}
\usepackage{amssymb}
\usepackage{tabularx}
\usepackage{enumerate}
\usepackage{graphicx}
\usepackage{texdraw}
\usepackage{color}
\usepackage{mathrsfs}
\usepackage{amsfonts,amssymb,amsmath}
\usepackage{epsfig}
\usepackage{adjustbox}

\makeatletter
\@addtoreset{equation}{section}

\makeatother
\newtheorem{thm}{Theorem}[section]
\newtheorem{lem}[thm]{Lemma}
\newtheorem{cor}[thm]{Corollary}
\newtheorem{prop}[thm]{Proposition}
\newtheorem{remark}[thm]{Remark}

\newcommand{\R}{\mathbb{R}}
\newcommand{\ve}{\varepsilon}

\begin{document}

\title[Mean Curvature Flow with Robin Boundary Conditions]
{Dimension-Dependent Asymptotic Dynamics for Mean Curvature Flow with Robin Boundary Conditions$^*$}
\thanks{$^*$ This research was partly supported by National Natural Science Foundation of China (No. 12471199, 11871148, 12001375).}
\author[X. Chen, B. Lou, X. Wang, L. Yuan]{Xinfu Chen$^\dag$, Bendong Lou$^{\ddag, \S}$, Xiaoliu Wang$^{\sharp}$ and Lixia Yuan$^{\ddag}$}
\thanks{$^\dag$ School of Mathematics, Southwestern University of Finance and Economics, Chengdu 611130, China.}
\thanks{$\ddag$ Mathematics and Science College, Shanghai Normal University, Shanghai 200234, China.}
\thanks{$\sharp$ School of Mathematics, Southeast University, Nanjing 210018, China.}
\thanks{$\S$ The corresponding author.}
\thanks{{\bf Emails:} {\sf xinfu@pitt.edu (X. Chen)}, {\sf lou@shnu.edu.cn (B. Lou)}, {\sf xlwang@seu.edu.cn (X. Wang)}, {\sf yuanlixia@shnu.edu.cn (L. Yuan)}}
\date{}


\maketitle

\begin{abstract}
We consider a graphical mean curvature flow in a cylinder with Robin boundary conditions, which arises as a geometric model for interface motion in the singular limit of the Allen--Cahn equation with nonlinear boundary conditions. It was shown in \cite{LWY} that, in the planar case, every solution converges to a translating Grim Reaper with a \emph{fixed profile} and \emph{finite speed}. In this paper, we investigate the radially symmetric problem in higher dimensions and reveal a completely different asymptotic dynamics caused by the spatial dimension. In contrast to the planar case, there is no fixed translating profile governing the long-time behaviour. Instead, the solution propagates with an exponentially increasing speed, while both the gradient $|Du|$ (away from the center) and the instantaneous speed $u_t$ diverge exponentially as $t\to\infty$. 
This reveals a fundamentally different asymptotic behaviour  induced by the interaction between the Robin boundary condition and the spatial dimension, that is, the translating profile continuously degenerates and becomes asymptotically ray-like. Since the equation becomes asymptotically degenerate and no uniform-in-time $C^0$, $C^1$, or $C^2$ estimates are available, our analysis relies on a new approach based on the zero number argument.
\end{abstract}

\baselineskip 16pt

\section{Introduction}

Given an embedded  $n$-dimensional hypersurface
$X_0:M^n\to\mathbb R^{n+1}$ in Euclidean space, we consider the one-parameter
family of hypersurfaces
$X(t):M^n\to\mathbb R^{n+1}$, $0<t<T$, generated by the mean curvature flow
(MCF), which is governed by
\[
\frac{\partial X(p,t)}{\partial t}
=\vec H,
\qquad
p\in M^n,\ \ 0<t<T.
\]
Geometrically, the mean curvature flow deforms a hypersurface in the direction
of its mean curvature vector $\vec H$, starting from the initial hypersurface
$X(0)=X_0$.

When the initial hypersurface is compact without boundary, the mean curvature
flow has been extensively studied; see, for example,
Gage and Hamilton \cite{GH1},
Grayson \cite{Grayson},
Huisken \cite{Huisken},
and many subsequent works.
The noncompact case has also attracted considerable attention; see
Chou and Zhu \cite{CZ1},
Ecker and Huisken \cite{EH1,EH2},
Cheng and Sesum \cite{Cheng1,Cheng-Sesum},
among others.
More recently,
Isenberg, Wu and Zhang \cite{IWZ2}
constructed a class of complete rotationally symmetric entire graphs with
prescribed growth at spatial infinity and showed that the corresponding flow
exists globally, escapes to spatial infinity as $t\to\infty$, while the second
fundamental form blows up at the rate
$(2t+1)^{(\gamma-1)/2}$ for each given $\gamma >0$.
In their earlier works \cite{IW,IWZ1},
they considered hypersurfaces contained in a cylinder and proved that the
corresponding flow develops finite-time singularities with the second
fundamental form blowing up at the rate $(T-t)^{-1}$.
For a comprehensive account of the mean curvature flow,
we refer the readers to the recent monograph by Andrews \cite{ACGL}.

These developments indicate that the long-time behaviour  of mean curvature flow
depends strongly on the geometry of the hypersurface and the boundary
conditions. In recent years, the asymptotic dynamics of geometric evolution
equations has become an active research topic, where translating solutions are
often found to determine the long-time behaviour  of the flow.
For the curve shortening flow,
Choi, Choi and Daskalopoulos \cite{CCD1}
proved the convergence of solutions with two ends asymptotic to parallel lines
to the translating Grim Reaper.
Later, they established an analogous convergence result for the Gauss
curvature flow \cite{CCD2}.
For a class of fully nonlinear one-dimensional parabolic equations,
including the curve shortening flow,
Kagaya and Liu \cite{KL1}
considered singular Neumann boundary conditions
($|Du|=\infty$ on the boundary)
and proved convergence to the corresponding translating solutions.
These results suggest that translating solutions often serve as asymptotic
attractors for geometric flows.

The above observations naturally raise the following question.
Does the same asymptotic picture remain valid for graphical mean curvature
flow with Robin boundary conditions?
For the planar problem the answer is affirmative, as was proved in
\cite{LWY}.
The main purpose of the present paper is to show that the higher-dimensional
problem exhibits a completely different asymptotic dynamics.
Indeed, the spatial dimension fundamentally changes the long-time behaviour  of
the flow, leading to a new asymptotic mechanism which is absent in the planar
case.

\vskip 10pt

We now formulate the problem.
Under the graphical mean curvature flow,
the graph function $u(x,t)$ satisfies
\begin{equation}\label{mcf1}
u_t=
\left(
\delta_{ij}
-
\frac{D_i uD_j u}{1+|Du|^2}
\right)
D_{ij}u,
\qquad
x\in\Omega,\ t>0,
\end{equation}
where $\Omega\subset\mathbb R^N$ is a bounded domain with smooth boundary.

When $N=1$, equation \eqref{mcf1} becomes
\begin{equation}\label{mcf2}
u_t=\frac{u_{xx}}{1+u_x^2},
\qquad
x\in(-1,1),\ t>0.
\end{equation}

In 1993,
Altschuler and Wu \cite{AW1}
studied translating solutions of \eqref{mcf2} under the boundary condition
\begin{equation}\label{b1}
u_x(\pm1,t)=\pm h,
\qquad
t>0,
\end{equation}
where $h>0$ is a constant.
This problem and its higher-dimensional analogue,
namely equation \eqref{mcf1} with the boundary condition
\eqref{b3} below,
have been extensively studied from the viewpoint of geometric analysis;
see
\cite{AW1,AW2,CW,Huisken1989,Guan,Mali,MWW}
and the references therein.
On the other hand,
these models also arise naturally as geometric descriptions of interface
motion in the singular limit of the Allen--Cahn equation with nonlinear
boundary conditions
(see
\cite{Matano2008,Chen1992,Fei,NMHS},
and the Appendix for a formal derivation of
\eqref{mcf1}, \eqref{b1} and \eqref{b3}).

Altschuler and Wu \cite{AW1}
proved that every solution of \eqref{mcf2}--\eqref{b1}
converges to the translating Grim Reaper
\begin{equation}\label{tw1}
\phi^c(x)+c(h)t
:=
-\frac1{c(h)}
\ln\!\bigl[ \cos(c(h)x)\bigr]
+c(h)t,
\end{equation}
where
$c(h):=\arctan h$.

A different phenomenon occurs when the boundary condition is replaced by the
Robin condition
\begin{equation}\label{b2}
u_x(1,t)=u(1,t),
\qquad
u_x(-1,t)=-u(-1,t),
\qquad
t>0.
\end{equation}
Chou and Wang \cite{CW}
proved that solutions may diverge to $\pm\infty$ as $t\to\infty$.
Later,
Lou, Wang and Yuan \cite{LWY}
showed that every global solution still converges locally to a translating
Grim Reaper, namely,
\begin{equation}\label{to-GR}
u(x,t+s)-u(0,s)
\to
-\frac2\pi
\ln\!\left(\cos\frac{\pi x}{2}\right)
+\frac\pi2 t
\qquad
\text{as }s\to\infty,
\end{equation}
in the $C^{2,1}_{\rm loc}((-1,1)\times\mathbb R)$ topology.
Unlike the problem with bounded boundary slope,
uniform-in-time gradient estimates are no longer available under the Robin
boundary condition.
Nevertheless,
the zero number argument developed in \cite{LWY}
yields uniform interior gradient estimates,
which are sufficient to recover the translating asymptotic profile.
This naturally raises the question whether the same asymptotic mechanism
persists in higher dimensions.
As we shall show, the answer is negative.

In higher dimensions, namely $\Omega\subset\R^N$ with $N\ge2$, the graphical
mean curvature flow is usually equipped with the capillary-type boundary
condition
\begin{equation}\label{b3}
  \frac{-Du\cdot\nu}{\sqrt{1+|Du|^2}}
  =g(x,t,u):=\cos\theta(x,t,u),
  \qquad
  x\in\partial\Omega,\ t>0,
\end{equation}
where $\nu$ denotes the inward unit normal vector to $\partial\Omega$, and
$\theta$ is the contact angle between the cylinder boundary
$\partial\Omega\times\R$ and the graph of $u(\cdot,t)$.
(See the Appendix for a formal derivation of
\eqref{mcf1} and \eqref{b3} from the singular limit of the Allen--Cahn equation
with nonlinear Robin boundary conditions.)

When $N=2$, $\Omega$ is strictly convex and $g=g(x)$ with sufficiently small
$|D_Tg|$, Altschuler and Wu \cite{AW2} proved that every solution converges
either to a minimal surface or to a translating solution.
More recently,
Ma, Wang and Wei \cite{MWW}
established uniform gradient estimates for general dimensions
under suitable geometric assumptions,
and consequently obtained convergence to translating solutions.
These results further support the general principle that,
under bounded boundary slopes,
the long-time dynamics is governed by a translating profile with fixed shape.

In the present paper we consider instead the Robin boundary condition
\begin{equation}\label{Robin-bdry-cond}
\frac{-Du\cdot\nu}{\sqrt{1+|Du|^2}}
=
\frac{u}{\sqrt{1+u^2}},
\qquad
x\in\partial\Omega,\ t>0,
\end{equation}
which is the natural higher-dimensional counterpart of
\eqref{b2}.
Moreover, we restrict ourselves to radially symmetric solutions.
More precisely,
we assume that $\Omega$ is the unit ball in $\R^N$ $(N\ge2)$ and
\[
u(x,t)=u(r,t),
\qquad
r=|x|\in[0,1].
\]
Then problem \eqref{mcf1}--\eqref{Robin-bdry-cond} reduces to
\begin{equation}\label{P}\tag{P}
\left\{
\begin{array}{ll}
\displaystyle
u_t
=
\frac{u_{rr}}{1+u_r^2}
+\frac{N-1}{r}u_r,
&
r\in(0,1),\ t>0,
\\[2mm]
u_r(0,t)=0,
\qquad
u_r(1,t)=u(1,t),
&
t>0.
\end{array}
\right.
\end{equation}
The corresponding mean curvature is
\[
H(r,t)
=
\frac{u_{rr}}{(1+u_r^2)^{3/2}}
+
\frac{(N-1)u_r}
{r\sqrt{1+u_r^2}}.
\]

Problem \eqref{P} is precisely the higher-dimensional analogue of the
one-dimensional problem \eqref{mcf2}--\eqref{b2} studied in \cite{LWY}.
In view of the planar convergence result, one may naturally expect that every
solution again approaches a translating solution with a fixed asymptotic
profile.
Surprisingly, this expectation is false.
Our analysis reveals a completely different asymptotic dynamics in higher
dimensions.
Unlike the planar case, there is no fixed translating profile governing the
long-time behaviour .
Instead, the profile continuously changes as the solution evolves,
the propagation speed increases exponentially,
and both the spatial gradient (away from the center) and the instantaneous
speed become unbounded as $t\to\infty$.
These phenomena show that the interaction between the Robin boundary condition
and the spatial dimension fundamentally changes the asymptotic dynamics of the
flow.
Our main result is the following.

\noindent
{\bf Main Theorem.}
{\it
Assume that $u_0\in C^2([0,1])$ is positive and satisfies the compatibility
conditions $u_0'(0)=0,\ u_0'(1)=u_0(1)$. Then problem \eqref{P} with initial data $u(r,0)=u_0(r)$ admits a unique global classical solution.
Moreover, as $t\to\infty$, the following assertions hold.

\medskip

\noindent
{\rm (i)}
$\displaystyle
\min_{r\in[0,1]}u(r,t)\to\infty,\qquad
\min_{r\in[0,1]}u_t(r,t)\to\infty,
$
and
\[
\min_{r\in[\varepsilon,1]}u_r(r,t)\to\infty
\qquad
\text{for every }
\varepsilon\in(0,1).
\]

\medskip

\noindent
{\rm (ii)}
\[
\frac{u_t(r,t)}{u(r,t)}
\stackrel{\rm w}{\longrightarrow}
N-1,
\qquad
rH(r,t)
\stackrel{\rm w}{\longrightarrow}
N-1
\]
in $C([0,1])$.

\medskip

In addition, when $N\geqslant 3$, the following sharper asymptotic estimates hold.

\medskip

\noindent
{\rm (iii)}
There exists a positive constant $C$ such that
\[
\left\|
\frac{u_r(r,t)}{u(r,t)}-r
\right\|_{C([0,1])}
\leqslant
Ce^{-\frac{N-1}{3}t},
\qquad
t\geqslant 0.
\]

\medskip

\noindent
{\rm (iv)}
There exist a constant $m_*\in\R$ and a function
$\vartheta(r,t)\to0$
uniformly on $[0,1]$ as $t\to\infty$ such that
\[
u(r,t)
=
e^{(N-1)t+\frac{r^2}{2}+m_*+\vartheta(r,t)},
\qquad
r\in[0,1],\ t\geqslant 0.
\]
}

\medskip
\begin{remark}\rm 
Unlike the planar problem studied in \cite{LWY}, the asymptotic behaviour  of
problem \eqref{P} is no longer governed by a translating solution with a fixed
profile.
Instead, both the propagation speed and the profile evolve continuously in
time.
In particular, not only the solution itself, but also its spatial gradient
(except at the center) and the instantaneous speed grow exponentially as
$t\to\infty$.
Consequently, there are no uniform-in-time
$C^0$, $C^1$ or $C^2$ estimates, and the corresponding asymptotic equations
become degenerate.
This reveals a fundamentally different asymptotic dynamics induced by the
interaction between the Robin boundary condition and the spatial dimension.
\end{remark}

\begin{remark}\rm
The lack of global gradient bounds also makes the classical gradient estimate
approach developed in
\cite{AW1,AW2,Huisken1989,Guan,Mali,MWW}
inapplicable.
Our analysis therefore relies on a different strategy based on the zero number
argument, which provides the crucial interior estimates required to determine
the asymptotic behaviour .
This approach may be useful for other geometric evolution equations with
unbounded boundary slopes.
\end{remark}

\begin{remark}\rm
The restriction $N\geqslant 3$ in assertions {\rm (iii)} and {\rm (iv)} is purely
technical.
It is only used in the proof of Lemma~\ref{lem:zeta-lower-est-1}.
We expect that this assumption can be removed; see
Remark~\ref{rem:N>2} for further discussion.
\end{remark}

The remainder of this paper is organized as follows.
Section~2 establishes the global existence together with several a priori
estimates.
Section~3 investigates the asymptotic behaviour  of $u$ and the normalized
gradient $u_r/u$.
Section~4 studies the asymptotics of $u_t$ and $u_t/u$.
Section~5 identifies the asymptotic limit of the mean curvature.
Finally, the Appendix presents a formal derivation of the graphical mean
curvature flow equation and the Robin boundary condition from the singular
limit of the Allen--Cahn equation with nonlinear boundary conditions.


\section{Global Existence}

This section establishes the global well-posedness of problem~\eqref{P}.
The main difficulty lies in the lack of uniform {\it a priori} estimates caused by the Robin boundary condition, under which neither the solution nor its gradient remains uniformly bounded as time evolves.
We first derive several fundamental estimates that capture the exponential growth of the solution while remaining sufficiently strong for the subsequent asymptotic analysis.
These estimates then allow us to prove the global existence of classical solutions.

\subsection{A priori estimates}

We begin with several basic estimates for the solution.
Although uniform-in-time $C^0$ and gradient bounds are unavailable,
the transformed variable introduced below enables us to obtain sharp exponential upper and lower bounds, which will play a central role throughout the paper.

Assume $u_0\in C^2([0,1])$, $u_0(r)>0$ and it satisfies the compatibility conditions. Assume further that the problem  \eqref{P} with initial data $u(r,0)=u_0(r)$ has a classical solution $u$ in  maximal existence interval $[0,T_\infty)$.

Motivated by the expected exponential growth rate obtained in the Main Theorem,
we introduce the normalized variable
$$
\displaystyle v(r,t) :=  \ln u - (N-1)t -\frac{r^2}{2} \quad \Leftrightarrow \quad
u(r,t)=e^{(N-1)t+\frac{r^2}{2}+v(r,t)}.
$$
Then $v$ solves

\begin{equation*}
\left\{
\begin{array}{ll}
 \displaystyle {\mathcal{N}_1}v:= v_t - \frac{1+v_{rr}+(r+v_r)^2}{1+u^2(r+v_r)^2} - \frac{N-1}{r}v_r=0, &  r\in(0,1),\ t\in (0,T_\infty),\\
v_r(0,t) = 0, \ \ v_r(1,t)=0, & t\in (0,T_\infty),\\
v(r,0)=v_0(r):= \ln u_0(r) - \frac{r^2}{2}, & r\in[0,1].
\end{array}
\right. \tag{Pv}
\end{equation*}

Define
\begin{equation}\label{def-m-M}
m(t) := \min_{r\in[0,1]}v(r,t),\quad
M(t) := \max\limits_{r\in [0,1]} v(r,t),\quad t\geqslant 0,
\end{equation}
and, for each $\tau\in [0,T_\infty)$, define
\begin{equation}\label{def of mt}
\underline{v}(r,t) :=m(\tau), \quad r\in[0,1],\ t \in [\tau,T_\infty) .
\end{equation}
One can verify that $\underline{v}$ is a subsolution in the time interval $[\tau, T_\infty)$. Hence,
$$
v(r,t)\geqslant  \underline{v}(r,t)=m(\tau), \quad r\in[0,1],\ t\in [\tau,T_\infty).
$$
This implies that the minimum $m(t)$ of $v(\cdot, t)$ satisfies
\begin{equation}\label{mono-mini-v}
m(0)\leqslant  m(t_1)\leqslant  m(t_2), \quad 0\leqslant  t_1\leqslant  t_2<T_\infty.
\end{equation}

The monotonicity of $m(t)$ immediately yields a lower exponential bound for the solution.

\begin{lem}\label{lem:lower-est-u}
The solution $u$ satisfies
$$
u(r,t)\geqslant  e^{(N-1)t+m(0)+\frac{r^2}{2}}, \quad r\in[0,1],\ t\in [0,T_\infty).
$$
\end{lem}

A complementary upper bound can be obtained by constructing a suitable supersolution, to say, $\bar{v}:=M(0)+ 2 t$ is a supersolution of the problem (Pv). Thus, we have the following upper estimate for $u$.

\begin{lem}\label{lem:upper-est-u}
The solution $u$ satisfies
$$
u(r,t)\leqslant  e^{(N+1) t + M(0) + \frac{r^2}{2}}, \quad r\in[0,1],\ t\in [0,T_\infty).
$$
\end{lem}

We next turn to gradient estimates.
Unlike the planar Robin problem studied in \cite{LWY}, global and even interior gradient bounds are expected to fail in higher dimensions.
The following computation nevertheless provides an evolution equation that will be repeatedly used in the sequel. More precisely, differentiating \eqref{mcf1} with respect to $x_k$ we obtain the equation about $D_k u$:
\begin{equation}
(D_k u)_t  =  \left( \delta_{ij} -\frac{D_i u D_j u}{1+|Du|^2}\right)  D_{ijk} u  + 2 \frac{D_i u D_j u Du \cdot DD_k u - D_{ik}u D_j u (1+|Du|^2)}{(1+|Du|^2)^2} D_{ij}u.
\end{equation}
Keeping in mind that $|D_k u| \leqslant  |u_r| = |u|$ on $\partial \Omega$, by the maximum principle, we obtain the following gradient estimate.
\begin{lem}
For any $T\in (0,T_\infty)$, there holds
$$
|D_k u (x,t)| \leqslant  \max\{\|D u_0\|_{L^\infty(\Omega)}, \|u\|_{L^\infty(\partial \Omega\times[0,T])}\},\quad x\in \Omega,\ t\in [0,T].
$$
\end{lem}

\subsection{Global existence of the solution}
To derive the global existence for a classical solution of a quasilinear parabolic equation, after obtaining the $L^\infty$ estimate and the gradient estimate one usually needs to establish the $C^\alpha$ estimate for $D u$, which can be attained if the equation satisfies some particular structure  and boundary conditions (cf. Lieberman \cite{Lie}).  For our problem, we refer to an existence result own to Ural'tseva \cite[Theorem 6]{Ur}, where a fully nonlinear parabolic type equation with a general boundary condition is considered. Since both of the equation and the boundary condition in our problem satisfy the structure conditions in \cite{Ur}, we can employ the result to obtain the following conclusion.

\begin{lem}
The problem \eqref{P} with initial data $u(r,0)=u_0(r)$ as stated above has a classical solution $u$ for all $t\geqslant  0$.
\end{lem}

The estimates established in this section guarantee the global existence of solutions and provide the basic quantitative information needed in the asymptotic analysis carried out in the subsequent sections.

\section{Asymptotic Behaviour I: The Asymptotic Profile and the Normalized Gradient}

The goal of this section is to identify the asymptotic profile of the solution.
As shown in the Main Theorem, the solution does not converge to a fixed translating profile as in the planar case.
Instead, after a suitable normalization, the solution exhibits an exponentially growing profile determined by the spatial dimension. 
Motivated by this exponential growth, we introduce a normalized quantity $\frac{u_r}{u}$. 
We first analyze this normalized gradient and then recover the asymptotic behaviour  of $u$ itself.

The next lemma provides the first control of the normalized gradient.

\begin{thm}\label{mainthm1}
Assume $N\geqslant 3$. Then the following conclusions hold for the global classical solution $u$.
\begin{enumerate}[{\rm (i).}]
\item
There exists a constant $C>0$ such that
$$
\Big|\frac{u_r}{u}-r\Big| \leqslant  Ce^{-\frac{N-1}{3}t}, \quad r\in[0,1],\ t\geqslant  0.
$$

\item There exists $m_*\in \R$ such that $\vartheta(r,t):=\ln u(r,t) - (N-1)t -\frac{r^2}{2} -m_*$ satisfies
$$
\lim_{t\to\infty}\sup_{0\leqslant  r \leqslant  1} \left|\vartheta(r,t) \right| = 0,
$$
that is, $\vartheta(r,t) = o(1)\ (t\to \infty)$ and 
$$
u(r,t) = e^{(N-1)t + \frac{r^2}{2} + m_* + \vartheta(r,t)},\qquad r\in [0,1], \ t\geq 0.
$$
\end{enumerate}
\end{thm}

\noindent
This theorem will be proved via several steps in the following subsections.
Before Lemma \ref{lem:zeta-upper-est} we only need $N\geqslant 2$. From Lemma \ref{lem:zeta-lower-est-1} to the end of this section, however, we require $N\geqslant 3$ for some technical reason (see Remark \ref{rem:N>2}).

\subsection{Monotonicity of $u(\cdot,t)$}
In this part we will use the zero number argument to show that  $u_r(r,t)$ will be positive (even larger) for large $t$, though $u'_0(r)$ is not necessarily to be so.

In order to construct suitable functions to compare with $u$, we first prepare some auxiliary functions. It was shown in Lou and Yuan \cite{LY} that, for any $k>0$, the equation in \eqref{P} has a translating solution $u= \varphi(r,k)+kt$ with
\begin{equation}\label{tw3-1}
\varphi (0,k)=\varphi_r (0,k)=0, \quad \varphi_{rr} (0,k)=\frac{k}{N},\quad \varphi_r(r,k)>0, \quad
\varphi_{rr}(r,k)>0 \ \text{for} \ r>0,
 \end{equation}
and
\begin{equation}\label{tw3-2}
  \varphi (r,k) =\left[ \frac{k}{2(N-1)} + o(1)\right] r^2  \quad\text{as} \quad r\to\infty.
 \end{equation}
In addition, $\varphi(r,k)$ satisfies
\begin{equation}\label{sol_tr}
\left\{\begin{array}{ll}
\displaystyle \varphi_{rr}(r,k) = \Big[k-\frac{N-1}{r}\varphi_r(r,k)\Big]\Big[1+\varphi_r^2(r,k)\Big], & r>0,\\
\displaystyle \varphi(r,k) = \Big[\frac{k}{2N}+o(1)\Big] r^2, & r\to 0.
\end{array}
\right.
\end{equation}

\begin{lem}\label{lem:prop-psi}
The function $\varphi(r,k)$ given above satisfies
$$
\frac{k}{N-1}r \geqslant  \varphi_r(r,k) \geqslant  \frac{k}{N}r,\quad r\in [0,1],
$$
and
$$
\frac{k}{2(N-1)}r^2 \geqslant  \varphi(r,k) \geqslant  \frac{k}{2N}r^2, \quad r\in[0,1].
$$
\end{lem}

\begin{proof}
The lemma can be proved easily by observing the slope fields on the lines $\varphi_r = \frac{k}{N}r$ and $\varphi_r = \frac{k}{N-1}r$.
\end{proof}

For any $k>0$ and any $p\geqslant  \varphi_r(1,k) -\varphi(1,k)$, set
\begin{equation}\label{def-ukp}
u_{k,p}(r,t) := \varphi(r,k)+kt+p.
\end{equation}
Then a direct verification shows that $u_{k,p}$ satisfies
\begin{equation}\label{p-ukm}
\left\{\begin{array}{ll}
\displaystyle (u_{k,p})_t = \frac{(u_{k,p})_{rr}}{1+((u_{k,p})_r)^2} + \frac{N-1}{r}(u_{k,p})_r, & 0<r<1,\ t >0,\\
(u_{k,p})_r(0,t) = 0, & t> 0,\\
(u_{k,p})_r(1,t) - u_{k,p}(1,t) = \varphi_r(1,k) -\varphi(1,k)-kt - p  <0, & t> 0.
 \end{array}
\right.
\end{equation}

\begin{lem}\label{lemma5}
Let $p(k):=\max \left\{\max\limits_{0\leqslant  r\leqslant 1}\{u_0(r)-\varphi(r,k)\}, \varphi_r(1,k)-\varphi(1,k)\right\}$. For every $p>p(k)$  and $t>0$, one of the following conclusions holds.
\begin{enumerate}[{\rm (i).}]
\item $u_{k,p}(\cdot,t) > u(\cdot,t)$  for $r\in[0,1]$;
\item $u_{k,p}(\cdot,t) < u(\cdot,t)$  for $r\in[0,1]$;
\item there exists a $z(t)\in [0,1]$ such that
\begin{equation}
\left\{\begin{array}{ll}
u_{k,p}(\cdot,t) > u(\cdot,t), & r\in[0,z(t)),\\
u_{k,p}(\cdot,t) < u(\cdot,t), & r\in(z(t),1],\\
u_{k,p}(z(t),t) = u(z(t),t).
 \end{array}
\right.
\end{equation}
\end{enumerate}
\end{lem}

\begin{proof}
Fix $p> p(k)$. Set $\eta := u-u_{k,p}$. Then $\eta(\cdot,0)<0$ on $[0,1]$. Set
$$
\begin{array}{l}
{\mathcal{A}} :=\{t\geqslant  0 \mid \eta(\cdot,t)<0 \ {\mathrm{on}} \ [0,1]\},\\
{\mathcal{B}} :=\{t\geqslant  0 \mid  \eta(\cdot,t)>0 \ {\mathrm{on}} \ [0,1]\},\\
{\mathcal{C}} := \left\{ t\geqslant 0 \left|
 \begin{array}{l}
 \mbox{there exists } z\in[0,1] \mbox{ such that } \eta(z,t)=0,\\
 \eta(\cdot,t)<0 \ {\mathrm{on}} \ [0,z),\  \eta(\cdot,t)>0 \ {\mathrm{in}} \ (z,1]
 \end{array}
  \right.
  \right\},
\end{array}
$$
and
$$
t_* := \sup\{t>0 \mid [0,t)\subset{\mathcal{A}}\cup{\mathcal{B}}\cup{\mathcal{C}}\}.
$$
Since ${\mathcal{A}}$ is open and $0\in {\mathcal{A}}$, we see that $t_*>0$. We claim
that $t_*=\infty$. Suppose by contradiction $t_*<\infty$. Consider three cases:
$$
(1)\  \eta(1,t_*)>0;\qquad \  (2)\ \eta(1,t_*)<0;\qquad \ (3)\ \eta(1,t_*)=0.
$$

In the case $\eta(1,t_*)>0$, there exists a small $\varepsilon>0$ such that $\eta(1,t)>0$ for $t\in[t_*-\varepsilon,t_*+\varepsilon]$. This implies that
$[t_*-\varepsilon,t_*)\subset {\mathcal{B}}\cup{\mathcal{C}}$. By the standard zero number diminishing properties (see, for example, Angenent \cite{Ang}), the zero number of $\eta(\cdot,t)$ does not increase in $t\in [t_*-\varepsilon,t_*+\varepsilon]$.
Hence, $[t_*-\varepsilon,t_*+\varepsilon]\subset {\mathcal{B}}\cup{\mathcal{C}}$.

In the case $\eta(1,t_*)<0$, we have $\eta(1,t)<0$ for $t\in[t_*-\varepsilon,t_*+\varepsilon]$ for some $\varepsilon>0$. As $[0,t_*)\subset{\mathcal{A}}\cup{\mathcal{B}}\cup{\mathcal{C}}$, we see that
$[t_*-\varepsilon,t_*)\subset{\mathcal{A}}$, and so $\eta(\cdot,t_*-\varepsilon)<0$ on $[0,1]$. Thus, $\eta<0$ on $[0,1]\times[t_*-\varepsilon,t_*+\varepsilon]$, which implies that
$[t_*-\varepsilon,t_*+\varepsilon]\subset{\mathcal{A}}$.

Finally, we consider the case $\eta(1,t_*)=0$. In this case $\eta_r (1,t_*)>0$ by the boundary conditions in \eqref{P} and \eqref{p-ukm}.  Hence, there exist a $\varepsilon\in(0,t_*)$ and a $\delta\in(0,1)$ such that
$$
\eta_r (r,t)>0 \mbox{ for }t\in[t_*-\varepsilon, t_*+\varepsilon] \mbox{ and } r\in[1-\delta,1],\qquad
\eta(1-\delta,t)<0 \mbox{ for } t\in[t_*-\varepsilon, t_*+\varepsilon].
$$
So we have $[t_*-\varepsilon, t_*)\subset{\mathcal{A}}\cup{\mathcal{C}}$. This implies that $\eta(r,t_*-\varepsilon)<0$ for any $r\in[0,1-\delta]$. By comparison we have
$\eta(r,t)<0$ for any $r\in[0,1-\delta]$ and $t\in[t_*-\varepsilon,t_*+\varepsilon]$. Since $\eta_r >0$ for $r\in[1-\delta,1]$ and $t\in[t_*-\varepsilon,t_*+\varepsilon]$. We see that
$[t_*-\varepsilon,t_*+\varepsilon]\subset{\mathcal{A}}\cup{\mathcal{C}}$.

In conclusion, $[t_*-\varepsilon,t_*+\varepsilon]\subset{\mathcal{A}}\cup{\mathcal{B}}\cup{\mathcal{C}}$.
This contradicts to the definition of $t_*$ and the assumption $t_*<\infty$. Consequently,  $t_*=\infty$ and ${\mathcal{A}}\cup{\mathcal{B}}\cup{\mathcal{C}}=[0,\infty)$. This proves the lemma.
\end{proof}

Now we show the monotonicity of $u(\cdot,t)$ for large $t$ by showing a positive lower bound for $u_r$.

\begin{lem}\label{lem:lower-est-ur}
There exists $T_0>0$ such that
$$
u_r(r,t) \geqslant  \frac{K(t)}{N}r, \quad r\in[0,1],\ t\geqslant T_0,
$$
where
\begin{equation}\label{def-K}
K(t):= \frac{e^{m(0)+(N-1)t}}{t+\frac{1}{N-1}}-1,\quad t\geqslant T_0.
\end{equation}
\end{lem}

\begin{proof}
Choose $T_0>0$ large such that
$$
K(t) = \frac{e^{m(0)+(N-1)t}}{t+\frac{1}{N-1}}-1 >2(N-1) \|u_0\|_{L^\infty},\quad t\geqslant T_0.
$$
For any fixed $t_0 \geqslant T_0$, we have
$$
2(N-1)\|u_0\|_{L^\infty}  <  k_0 := K(t_0) <  \frac{e^{m(0)+(N-1)t_0}}{t_0+\frac{1}{N-1}}.
$$
It follows that
$$
e^{m(0)+(N-1)t_0}-k_0 t_0 > \frac{k_0}{N-1} >\|u_0\|_{L^\infty}+\frac{k_0}{2(N-1)}.
$$
Using the properties of $\varphi$ in Lemma \ref{lem:prop-psi} we have
$$
\left\{\begin{array}{ll}
e^{m(0)+(N-1)t_0 }-k_0 t_0-\varphi(1,k_0) >  \varphi_r(1,k_0) - \varphi(1,k_0), \\
e^{m(0)+(N-1)t_0 }-k_0 t_0-\varphi(1,k_0) >  \|u_0\|_{L^\infty}.
\end{array}
\right.
$$
Consequently,
\begin{equation}\label{>p(t_0)}
e^{m(0)+(N-1)t_0} - k_0 t_0-\varphi(1,k_0) > p(k_0) :=  \max \Big\{\max\limits_{0\leqslant  r\leqslant 1}\{u_0(r)-\varphi(r,k_0)\}, \varphi_r(1,k_0)-\varphi(1,k_0)\Big\}.
\end{equation}

For any $r_0\in (0,1]$, set
$$
\tilde{p}(r_0, t_0) := u(r_0,t_0) - [k_0 t_0 +\varphi(r_0,k_0)].
$$
Then, \eqref{>p(t_0)} and $u(r_0,t_0)\geqslant  e^{m(0)+(N-1)t_0}$ in Lemma \ref{lem:lower-est-u} imply that $\tilde{p}(r_0, t_0)  > p(k_0)$. Define $u_{k_0,\tilde{p}}(r,t)$ as in \eqref{def-ukp}. Then we have
$$
u(r_0,t_0)=u_{k_0,\tilde{p}}(r_0,t_0).
$$
By Lemma \ref{lemma5}, we have
$$
u(r,t_0) < u_{k_0,\tilde{p}}(r,t_0) \mbox{ for } r\in[0,r_0),\quad u(r,t_0) > u_{k_0,\tilde{p}}(r,t_0) \mbox{ for } r\in(r_0,1],
$$
and so
\begin{equation}\label{grad_est}
u_r(r_0,t_0)\geqslant  \varphi_r(r_0,k_0)\geqslant  \frac{k_0}{N}r_0 = \frac{K(t_0)}{N}r_0.
\end{equation}
This proves the lemma.
\end{proof}

Combining this lemma with the estimates in Lemmas \ref{lem:lower-est-u} and \ref{lem:upper-est-u} we have the following result.
\begin{cor}\label{cor:ur-u-begin-est}
There exist $T_0>0,\ \varepsilon_0>0, \ K_0>0$  such that
$$
u(r,T_0+t) \geqslant  e^{(N-1)t}, \quad r\in[0,1],\ t\geqslant 0,
$$
and, at the moment $t=T_0$,
$$
\varepsilon_0 r \leqslant  \frac{u_r(r,T_0)}{u(r,T_0)} \leqslant  K_0 r,\quad r\in [0,1].
$$
\end{cor}
\noindent
Since we will study the asymptotic behaviour of $u$ in the rest of the paper, without loss of generality we shift the time such that the conclusions in this corollary hold for $T_0 = 0$.

\subsection{Finer upper bound of $v$}
Define $v$ by $u(r,t) = e^{(N-1)t+\frac{r^2}{2}+v(r,t)}$ as before. Then $v$ solves the problem (Pv) in
$t\in [0,\infty)$. Recall that we have obtained an upper bound $M(0)+2t$ for $v$ in Lemma \ref{lem:upper-est-u}. Since $M(0)+2t$ is unbounded, it is not good enough to give further properties for $v$. In this subsection
we will present a finer upper bound for $v$ by constructing a supersolution.

For any $k >0$ and $t_0\geqslant 0$, define
$$
q(r,k):= \frac{1}{N}\int_r^1\frac{\rho(e^{-k \rho}-e^{-k })}{1-e^{-k }} d\rho,
$$
and
$$
\bar{v}(r,t) := q(r,k) + L [e^{-(N-1)t_0}-e^{-(N-1)t}] + M(t_0),\quad r\in [0,1],\ t\geqslant t_0,
$$
where $M(t_0)$ is defined by \eqref{def-m-M} and $L>0$ is a constant to be determined below.
Then, with $\alpha := \frac{1}{N(1-e^{-k })}$, we have
$$
0\geqslant  \bar{v}_r = -\alpha r[e^{-k  r}-e^{-k  }] \geqslant  -\alpha r[1-e^{-k  }] = -\frac{r}{N},
$$
\begin{equation}\label{r+v_r}
r+\bar{v}_r = r[1-\alpha (e^{-k  r}-e^{-k  }) ] \geqslant  \frac{N-1}{N}r\geqslant 0,
\end{equation}
and $\bar{v}_r(0,t)=\bar{v}_r(1,t)=0$. Thus
\begin{equation*}
\begin{array}{lll}
\displaystyle {\mathcal{N}_1}\bar{v} & = & \displaystyle \frac{1}{1+u^2(r+\bar{v}_r)^2} (r+\bar{v}_r)^2\Big[(N-1) L e^{-(N-1)t}u^2-\frac{N-1}{r}\bar{v}_ru^2 -1\Big] \\
\displaystyle & & + \displaystyle \frac{1}{1+u^2(r+\bar{v}_r)^2} \left\{ (N-1) L e^{-(N-1)t}
-1-\bar{v}_{rr}-\frac{N-1}{r}\bar{v}_{r}\right\}.
\end{array}
\end{equation*}
Assume that $L\geqslant  2$. Then by Corollary \ref{cor:ur-u-begin-est} we have
$$
L e^{-(N-1)t}u^2 \geqslant  L e^{(N-1)t} \geqslant \frac{L}{2}e^{(N-1)t}+1.
$$
Combining with $\bar{v}_r\leqslant  0$ we have
\begin{equation*}
\begin{array}{lll}
{\mathcal{N}_1}\bar{v} & \geqslant & \displaystyle \frac{\frac{N-1}{2} L e^{(N-1)t}(r+\bar{v}_r)^2 +[(N-1) L e^{(N-1)t} -1-\bar{v}_{rr}-\frac{N-1}{r}\bar{v}_r]}
{1+u^2(r+\bar{v}_r)^2} \\
& \geqslant & \displaystyle \frac{\sqrt{2}(N-1) L (r+\bar{v}_r) - [1+\bar{v}_{rr}+\frac{N-1}{r}\bar{v}_r]} {1+u^2(r+\bar{v}_r)^2}.
\end{array}
\end{equation*}
By the definitions of $\bar{v},\ \alpha$ and by \eqref{r+v_r} we have
\begin{eqnarray*}
\displaystyle 1+\bar{v}_{rr}+\frac{N-1}{r}\bar{v}_r & = & 1 -\alpha (e^{-k  r}-e^{-k  }) + \alpha k  re^{-k  r} - (N-1)\alpha (e^{-k  r}-e^{-k  })\\
&= & \alpha N (1-e^{-k  r}) + \alpha k  re^{-k  r}\\
& \leqslant & (r+\bar{v}_r)\Big[\frac{N}{r}(1-e^{-k  r})+k  e^{-k  r} \Big]\cdot \frac{N\alpha}{N-1}\\
& \leqslant & (r+\bar{v}_r)\cdot \frac{N\alpha}{N-1}(N+1)k   \\
& = &\frac{(N+1)k  }{(N-1)(1-e^{-k})}(r+\bar{v}_r),\quad r\in [0,1],\ t\geqslant t_0.
\end{eqnarray*}
If we take
$$
L = L(k):= \max\Big\{\frac{(N+1)k  }{\sqrt{2}(N-1)^2 (1-e^{-k  })},\ 2\Big\},
$$
then we have
$$
[{1+u^2(r+\bar{v}_r)^2}]\cdot {\mathcal{N}_1 }\bar{v} \geqslant  (r+\bar{v}_r)\Big[\sqrt{2}(N-1)L - \frac{(N+1)k  } {(N-1)(1-e^{-k  })}\Big]\geqslant  0,\quad r\in [0,1],\ t\geqslant t_0.
$$
Therefore, $\bar{v}$ is a supersolution in $[t_0,\infty)$, and by comparison we prove the following lemma.

\begin{lem}\label{lem:fine-upper-v}
For each $k  >0$ and $t_0\geqslant  0$, there holds
$$
v(r,t) \leqslant q(r,k) + L (k )[e^{-(N-1)t_0}-e^{-(N-1)t}] + M(t_0), \quad r\in [0,1],\ t\geqslant  t_0.
$$
\end{lem}

\noindent
This lemma gives a finer upper bound for $v$. Taking $k=1$ and combining with \eqref{mono-mini-v} we have the following result.

\begin{cor}\label{cor:C0est-v-u}
There exists a constant $C>0$ depending only on $N$ such that
$$
m(0):=\min_{r\in[0,1]} v_0(r) \leqslant   v(r,t)  \leqslant  M(0)+C := \max\limits_{r\in [0,1]} v_0(r) + C,\quad r\in [0,1],\ t\geqslant 0.
$$
Consequently,
$$
e^{(N-1)t+\frac{r^2}{2}+m(0)}\leqslant u(r,t)\leqslant e^{(N-1)t+\frac{r^2}{2}+M(0)+C},\quad r\in [0,1],\ t\geqslant 0.
$$
\end{cor}

\begin{cor}\label{conv_v}
For $m(t)$ and $M(t)$ defined in \eqref{def-m-M}, there exists $m_*$ and $M_*$ such that
$$
\lim_{t\to\infty}m(t)=m_*, \quad \lim_{t\to\infty} M(t)=M_*.
$$
\end{cor}

\begin{proof}
Since $m(t)$ is an increasing function and it is bounded, we see that there exists $m_*$ such that $\lim\limits_{t\to\infty} m(t)=m_*$.

Next, set $M_*:=\liminf\limits_{t\to\infty} M(t)$. Fix any small $\varepsilon>0$, then there exists a time sequence $\{t_j\}_{j=1}^\infty\subset[0,\infty)$ such that
$$
M(t_j) \leqslant  M_*+\varepsilon.
$$
Note that $\lim\limits_{k\to\infty}\|q(\cdot,k)\|_{L^\infty}=0.$ There exists a large $k= k_\varepsilon>0$ such that $\|q(\cdot,k_\varepsilon)\|_{L^\infty}\leqslant  \varepsilon$. Using the estimate in Lemma \ref{lem:fine-upper-v} we have
$$
M(t) \leqslant \|q(\cdot, k_\ve)\|_{L^\infty} + L (k_\ve) e^{-(N-1)t_j} + M(t_j) , \quad t\geqslant  t_j.
$$
Thus,
$$
\limsup\limits_{t\to\infty} M(t) \leqslant  M_* + 2\ve + L(k_\ve)e^{-(N-1)t_j}.
$$
By sending $j\to\infty$, we obtain the assertion. This proves the corollary.
\end{proof}

\subsection{Gradient estimate of $w$}
Define $w$ by $u=e^{(N-1)t+w}$. Then
$$
u_r=uw_r, \quad u_{rr}=u(w_{rr}+w_r^2).
$$
Actually, recalling the definition of $v$ we see that $w=v+\frac{r^2}{2}$.
A direct computation shows that
\begin{equation}\label{P-w}
 \left\{
 \begin{array}{ll}
 w_t = \displaystyle \frac{w_{rr}+w_r^2}{1+u^2w_r^2}+\frac{N-1}{r}w_r-(N-1), & r\in (0,1),\ t>0,\\
 w_r(0,t)=0,\quad w_r(1,t)=1, & t>0.
 \end{array}
 \right.
\end{equation}
Furthermore, if we set $\zeta :=w_r$, then $\zeta$ satisfies
\begin{equation}\label{P-zeta}
\left\{
\begin{array}{ll}
{\mathcal{N}}_2\zeta := \displaystyle \zeta_t - \frac{\zeta_{rr}}{1+u^2\zeta^2} - \frac{N-1}{r}\Big[\zeta_r-\frac{\zeta}{r}\Big] & \\
\qquad \ \ \ \  \displaystyle  +\frac{2\zeta}{(1+u^2\zeta^2)^2}\Big[
u^2(\zeta_r^2+\zeta^2\zeta_r+\zeta^4)-\zeta_r\Big]=0, & r\in (0,1),\ t>0,\\
\zeta(0,t)=0, \quad  \zeta(1,t)=1, & t>0.
\end{array}
\right.
\end{equation}

By Corollary \ref{cor:ur-u-begin-est} we assume without loss of generality that
\begin{equation}\label{assump-u0-zeta0}
u(r,t)\geqslant  e^{(N-1)t},\quad \varepsilon_0 r \leqslant  \zeta(r,0) \leqslant  K_0 r, \quad r\in[0,1],\ t\geqslant  0,
\end{equation}
for $K_0\geqslant 1 >\varepsilon_0 >0$.

We now give finer upper bound for $\zeta$.

(1) Set $\overline{\zeta}_1 (r,t) := K_0 r$. Then
$$
{\mathcal{N}}_2\overline{\zeta}_1  = \frac{2\overline{\zeta}_1 }{\Big(1+u^2\overline{\zeta}_1^2\Big)^2}
[u^2(K_0^2 +K_0^3 r^2+K_0^4 r^4)-K_0]\geqslant  0,\quad r\in [0,1].
$$
Thus ${\zeta}(r,t) \leqslant  \overline{\zeta}_1 (r,t) =K_0 r$ for all $r\in[0,1],\ t\geqslant  0.$

(2) Set $\overline{\zeta}_2 (r,t) := r+(K_0-1)e^{-(N-1)t}$. Then
$$
 {\mathcal{N}}_2\overline{\zeta}_2 = (N-1)(K_0-1)e^{-(N-1)t}\Big[-1+\frac{1}{r^2}\Big] + \frac{2\overline{\zeta}_2 [u^2-1+u^2(\overline{\zeta}_2^2+\overline{\zeta}_2^4)]}{\Big(1+u^2 \overline{\zeta}_2^2 \Big)^2}\geqslant  0.
$$
Thus $\zeta(r,t) \leqslant  r+(K_0-1)e^{-(N-1)t}$  for all $r\in[0,1],\ t\geqslant  0$.

Combining these two cases together, we obtain a finer upper bound for $\zeta$:
\begin{lem}\label{lem:zeta-upper-est}
Assume \eqref{assump-u0-zeta0} holds. Then
$$
\zeta(r,t) \leqslant  \min\left\{K_0 r,\ r+(K_0-1)e^{-(N-1)t}\right\}, \quad r\in[0,1],\ t\geqslant  0.
$$
\end{lem}

Next we give finer lower bound for $\zeta$.  For some positive $\varepsilon\in (0,1)$ and $\sigma\in (0,1)$ to be determined below, we set $\underline{\zeta}_1 := \varepsilon r^{1+\sigma}$.
Using the Cauchy inequality $1+(u\underline{\zeta}_1)^2\geqslant  2u\underline{\zeta}_1$, one can show that
\begin{eqnarray*}
{\mathcal{N}}_2\underline{\zeta}_1 & = & \displaystyle
-\frac{\varepsilon(1+\sigma)\sigma r^{\sigma-1}}{1+u^2\underline{\zeta}_1^2}
- (N-1)\varepsilon\sigma r^{\sigma-1} + \frac{2\underline{\zeta}_1}{[1+u^2 \underline{\zeta}_1^2]^2} \left[ u^2\big(\underline{\zeta}_{1r}^2 + \underline{\zeta}_1^2 \underline{\zeta}_{1r} + \underline{\zeta}_1^4 \big) - \underline{\zeta}_{1r} \right] \\
& \leqslant & \displaystyle   - \frac{\varepsilon(1+\sigma)\sigma r^{\sigma-1} u^2 \underline{\zeta}_1^2}{[1+u^2 \underline{\zeta}_1^2]^2}  -(N-1)\varepsilon\sigma r^{\sigma-1}  \\
& & + \frac{2u^2 \underline{\zeta}_1}{[1+u^2 \underline{\zeta}_1^2]^2} \left[ \ve^2 (1+\sigma)^2 r^{2\sigma} + \ve^3 (1+\sigma)r^{2+3\sigma} + \ve^4 r^{4+4\sigma} \right]\\
& = &  \displaystyle  -(N-1)\varepsilon\sigma r^{\sigma-1} \\
& &   +\frac{u^2 \underline{\zeta}_1^2}{[1+u^2 \underline{\zeta}_1^2]^2} \left[ -\ve \sigma(1+\sigma)r^{\sigma -1} + 2\ve (1+\sigma)^2 r^{\sigma -1} +2\ve^2 (1+\sigma) r^{1+2\sigma} +2\ve^3 r^{3+3\sigma}\right]\\
& \leqslant &  \displaystyle  -(N-1)\varepsilon\sigma r^{\sigma-1}  +\frac{u^2 \underline{\zeta}_1^2}{[1+u^2 \underline{\zeta}_1^2]^2} \left[ (\sigma^2 + 3\sigma +2) + 2\ve (1+\sigma) +2\ve^2\right] \ve r^{\sigma -1} \\
& \leqslant &  \displaystyle  -(N-1)\varepsilon\sigma r^{\sigma-1}  +\frac14 \left[ \sigma^2 + 3\sigma +2 + 6\ve \right] \ve r^{\sigma -1} \\
& =  & \frac{\ve r^{\sigma -1}}{4} H(N,\sigma,\ve) := \frac{\ve r^{\sigma -1}}{4} \left[\sigma^2 + (7-4N)\sigma +2 + 6\ve \right].
\end{eqnarray*}
When $N\geqslant 3$, we set
\begin{equation}\label{sigma-ve-1}
\sigma_0 := \frac{4N-7-\sqrt{16N^2 -56 N+41}}{2}, \qquad
\ve^0 := \min\left\{ \frac12,\ \frac{(4N-7)\sigma -\sigma^2 -2}{7} \right\}.
\end{equation}
Now if $N\geqslant 3$ and if we choose
\begin{equation}\label{sigma-ve-2}
\sigma_0 <\sigma <1\qquad \mbox{and}\qquad 0<\ve < \ve^0
\end{equation}
(which includes the special case $\sigma =\frac12$), then we have
${\mathcal{N}}_2\underline{\zeta}_1\leqslant 0$. By comparison we obtain the following lower bound for $\zeta$.

\begin{lem}\label{lem:zeta-lower-est-1}
Assume $N\geqslant 3$ and \eqref{assump-u0-zeta0} holds. Then for any  $\sigma\in (\sigma_0,1)$ and $\ve\in (0,\ve^0]$ for $\sigma_0,\ \ve^0$ defined by \eqref{sigma-ve-1}, there holds
$$
\zeta(r,t) \geqslant  \varepsilon r^{1+\sigma}, \quad
 r\in[0,1],\ t\geqslant  0.
$$
\end{lem}

Combining this lemma with the previous one we have the following corollary.
\begin{cor}\label{cor:u_r/ru}
Assume $N\geqslant 3$ and \eqref{assump-u0-zeta0} holds. Then
$$
\varepsilon_0 r^\sigma \leqslant \frac{\zeta}{r} = \frac{u_r}{ru}\leqslant K_0,\quad r\in (0,1],\ t\geqslant 0,
$$
where $\sigma\in (\sigma_0, 1)$ and $K_0\geqslant 1>\ve^0 >0$.
\end{cor}

\begin{remark}\label{rem:N>2}\rm
To prove Lemma \ref{lem:zeta-lower-est-1} we use an additional condition $N\geqslant 3$, which ensure that $H(N,\sigma,\ve)<0$  for some $\sigma\in (0,1)$ and $\ve\ll 1$. Note that this is the only place where we need this condition. We guess that this additional condition can be omitted.
\end{remark}

Based on Lemma \ref{lem:zeta-lower-est-1}, we can even give a better lower bound for $\zeta$. Define
$$
{\mathcal{N}}_3\phi := \phi_t - \frac{\phi_{rr}}{1+u^2\zeta^2}
- \frac{N-1}{r}\Big(\phi_r - \frac{\phi}{r}\Big)
+ \frac{2\zeta u^2}{(1+u^2\zeta^2)^2}[\phi_r^2 + \phi^2 \phi_r +\phi^4] - \frac{2\zeta\phi_r}{(1+u^2\zeta^2)^2}.
$$
For $t_0 > 0$ to be determined, set
$$
\underline{\zeta}_2 := r- e^{- \frac{(N-1)(t-t_0)}{3}}, \quad
Q :=\left\{ (r,t)\; \left|\; t>t_0,\ e^{- \frac{(N-1)(t-t_0)}{3}} < r \leqslant 1 \right. \right\}.
$$
Note that, for any $0<\theta<1$, by Young's inequality there holds
$$
(1+u^2\zeta^2)^2 \geqslant  2u\zeta (1+u^2\zeta^2) \geqslant  2 (u\zeta)^{2+\theta}.
$$
So, in the set $Q$,  by $u\geqslant e^{(N-1)t}$ and Lemma \ref{lem:zeta-lower-est-1} we have
\begin{eqnarray*}
{\mathcal{N}}_3\underline{\zeta}_2
&\leqslant & \displaystyle
e^{- \frac{(N-1)(t-t_0)}{3}} \Big[\frac{N-1}{3} -\frac{N-1}{r^2}\Big] + \frac{\zeta u^2}{(u\zeta)^{2+\theta}} (1+r^2+r^4)\\
&\leqslant  & \displaystyle  e^{-\frac{(N-1)(t-t_0)}{3}} \Big[\frac{N-1}{3} -\frac{N-1}{r^2}\Big] + \frac{3}{u^\theta \zeta^{1+\theta}}\\
&\leqslant  & \displaystyle  e^{-\frac{(N-1)(t-t_0)}{3}} \Big[\frac{N-1}{3}-\frac{N-1}{r^2}\Big] +
\frac{3 e^{-\theta(N-1)t}}{\varepsilon^{1+\theta} r^{(1+\sigma)(1+\theta)}}.
\end{eqnarray*}
If we take $\sigma=\frac12,\ \theta =\frac13$ and $t_0$ satisfying
$ 9 e^{-\frac{(N-1)t_0}{3}} = 2 \ve^{\frac43} (N-1)$, then
\begin{eqnarray*}
{\mathcal{N}}_3\underline{\zeta}_2
&\leqslant  & \displaystyle  e^{-\frac{(N-1)(t-t_0)}{3}} \left[\frac{N-1}{3} -\frac{N-1}{r^2} + \frac{3e^{-\frac{(N-1)t_0}{3}}}{\varepsilon^{\frac43 } r^2}\right]\\
& \leqslant & (N-1) e^{-\frac{(N-1)(t-t_0)}{3}} \frac{r^2 -1}{3r^2}\leqslant 0,\quad (x,t)\in Q.
\end{eqnarray*}
Using the comparison principle in $Q$ we derive the following result.

\begin{lem}\label{lem:zeta-lower-est}
Assume $N\geqslant 3$ and \eqref{assump-u0-zeta0} holds. Then there exists $t_0>0$ such that
$$
\zeta(r,t) \geqslant  r - e^{-\frac{N-1}{3}(t-t_0)}, \quad r\in[0,1],\ t>t_0.
$$
\end{lem}

\subsection{The proof of Theorem \ref{mainthm1}}
Based on the estimates in the previous subsections we prove Theorem \ref{mainthm1} in this part.

The conclusion (i) follows from Lemma \ref{lem:zeta-upper-est} and Lemma \ref{lem:zeta-lower-est}.

To show the conclusion (ii), we could employ Corollary \ref{conv_v} and to show $m_*= M_*$.
Recall that $m(t):= \min_{r\in[0,1]}v(r,t)$ is increasing and tends to $m_*$ as $t\to \infty$. So, for any $\varepsilon>0$, there exists a $T>0$ such that
$$
m(t) \leqslant  m_* + \varepsilon, \quad  t\geqslant  T.
$$
When $N\geqslant 3$, by Lemma \ref{lem:zeta-upper-est} and Lemma \ref{lem:zeta-lower-est} we have
$$
|v_r(r,t)| = |\zeta(r,t)-r| \leqslant Ce^{-\frac{N-1}{3}t} \leqslant \ve, \quad r\in[0,1],\ t\geqslant  T.
$$
Hence,
$$
v(r,t)\leqslant  m(t)+\ve \leqslant m_*+2\ve, \quad r\in [0,1],\ t\geqslant  T.
$$
Consequently,
$$
M_*\leqslant  m_*+2\ve.
$$
Sending $\ve\to 0$, we obtain $M_*=m_*$.

This completes the proof of Theorem \ref{mainthm1}. 
\qed

\medskip
The analysis in this section identifies the asymptotic profile of the solution and shows that the normalized gradient converges to the spatial variable $r$.
These results will be the key input for the study of the asymptotic propagation speed in the next section.

\section{Asymptotic Behaviour II: The Asymptotic Propagation Speed}

The purpose of this section is to determine the asymptotic behaviour  of the
instantaneous speed $u_t$. 

The results in the previous section identify the asymptotic profile of the
solution and show that the normalized gradient $\frac{u_r}{u}$ converges to
$r$. We now use these properties to analyze the propagation speed of the flow.

A remarkable feature is that, unlike the planar case where the translating
speed remains finite, the speed in higher dimensions grows exponentially in
time.

\subsection{The infinity limit of the instantaneous velocity}

\begin{lem}
There exists a time $T$, such that $u_t(r,t)>0$ for $r\in[0,1], t\geqslant  T$.
\end{lem}

\begin{proof}
We first show that $u_t(0,t)\to\infty$ as $t\to\infty$. For any given $k>0$, the translating solution $\varphi(r,k)+k t$ solves the equation in \eqref{P}, though it does not necessarily satisfy the Robin boundary condition. Then we can use the zero number argument to consider the number of the intersection points between $u(r,t)$ and $\varphi(r,k)+k t+h$ ($h>0$ is large), as it was shown in the proof of Lemma
\ref{lemma5} (see also \cite{LWY} Section 4), to conclude that
$$
u_r(r,t) \geqslant \varphi_r (r,k),  \qquad r\in (0,1],\ t\gg 1.
$$
So, for any large $t$, there exists $h(t)$ such that
$$
 \begin{array}{l}
u(r,t) > \varphi(r,k) + k t + h(t), \quad  r\in(0,1],\\
u(0,t) = \varphi(0,k) + k t + h(t),\\
u_r(0,t) = 0 = \varphi_r(0,k),
\end{array}
$$
that is, $\varphi(r,k) + k t + h(t)$ lies below $u(\cdot,t)$ and is tangent to $u(\cdot,t)$ at $r=0$. Hence their mean curvatures at $r=0$ satisfy
$$
\left. \frac{u_{rr}}{(1+u_r^2)^{3/2}}+\frac{(N-1)u_r}{r\sqrt{1+u_r^2}}\right|_{r=0}
\geqslant  \frac{\varphi_{rr}}{(1+\varphi^2_r)^{3/2}} + \left. \frac{(N-1)\varphi_r}{r\sqrt{1+\varphi^2_r}}\right|_{r=0}.
$$
This implies that
$$
u_{rr}(0,t) \geqslant  \varphi_{rr}(0,k).
$$
This is true for all large $t$, say, $t\geqslant  T(k)$. Then, using the equation of $u$ we have
$$
u_t(0,t) = Nu_{rr}(0,t)\geqslant  N\varphi_{rr}(0,k)=k, \quad t\geqslant  T(k).
$$
Taking $k$ to be larger and larger we conclude that
$$
u_t(0,t)\to\infty \mbox{\ \ as\ \ } \ t\to\infty.
$$

Next we show $u_t(r,t)>0$ for $r\in[0,1],\ t\geqslant  T$. Set $\xi :=u_t$. Then
\begin{equation}\label{uteqn}
\left\{
\begin{array}{ll}
\xi_t = \displaystyle \frac{\xi_{rr}}{1+u_r^2} + \Big[\frac{-2u_ru_{rr}}{(1+u_r^2)^2}+\frac{N-1}{r}\Big]\xi_r, &  r\in(0,1),\ t>0,\\
\xi_r(0,t)=0, \quad \xi_r(1,t)=\xi(1,t), & t>0.
\end{array}
\right.
\end{equation}
Note that if $\xi(r,t_0)=u_t(r,t_0)\geqslant, \not\equiv  0$ for some $t_0\geqslant  0$, then by the maximum principle we have
$$
\xi(r,t) >  0, \quad r\in[0,1],\ t>  t_0.
$$
This is the desired conclusion. We assume by contradiction that, for any $t>0$, there exists $y(t)\in[0,1]$ such that \begin{equation}\label{zeropt}
\xi(y(t),t) <  0.
\end{equation}
From the previous step, we know that $y(t)>0$ rather than $y(t)=0$ for all large time $t$, say, $t\geqslant  T_1$.
Set $Q_1:= [0,1]\times [T_1,\infty)$. Then this domain contains a connected component $Q_2$, including
the left boundary $\{0\}\times [T_1,\infty)$ of $Q_1$ such that $\xi(r,t)>0$ in $Q_2$.
Clearly, $Q_2 \subset,\ \not= Q_1$ by (\ref{zeropt}), and $\xi=0$ on the right boundary of $Q_2$. Now, for any large $T_2$, we use the maximum principle for $\xi$ in $Q_3:= Q_2\cap ([0,1]\times [T_1,T_2])$ to conclude that the maximum of $\xi$ in $Q_3$ is attained on the bottom $\overline{Q_2} \cap ([0,1]\times \{T_1\})$ or on the left boundary $\{0\}\times [T_1,T_2]$.
When $T_2$ is sufficiently large, we have
$$
\max\limits_{T_1\leqslant  t\leqslant  T_2} \xi(0,t) > \max\limits_{0\leqslant  r\leqslant  1}\xi(r,T_1)
$$
due to $\xi(0,t)\to\infty\ (t\to \infty)$ in the previous step. Hence, the (positive) maximum of $\xi$ in $Q_3$ is attained at $(0, \tau )$ for some $T_1 < \tau \leqslant  T_2$. Using the Hopf Lemma at this point have
$$
\xi_r(0,\tau)<0,
$$
which contradicts the boundary condition $\xi_r(0,t)\equiv 0$.
This proves the lemma.
\end{proof}

\begin{thm}\label{ut_conv}
As $t\to \infty$, there holds $s(t):=\min\limits_{r\in [0,1]}u_t(r,t)\to\infty$.
\end{thm}

\begin{proof}
Using the maximum principle for $\xi$, it is easily seen that, when $t$ is large, $s(t)$ is positive and is increasing in $t$. Suppose by contradiction that, for some $S^*\in (0,\infty)$,
$$
s(t) < S^*, \quad  t>0.
$$
This implies that, for all $t>0$, there exists $z(t)\in[0,1]$ such that
$$
\xi (z(t),t) < S^* .
$$
Without loss of generality, assume that $z(t)$ is the smallest one of such points in $[0,1]$.
Since $\xi(0,t)\to\infty\ (t\to\infty)$, we have $\xi(0,t)>S^*$ for large $t$, say
$t\geqslant T_3$. Then $z(t)\in (0,1]$ for $t\geqslant T_3$. Using a similar argument as in the proof of the previous lemma in the domain
$$
Q_4 := \{(r,t) \mid 0<r<z(t), \ t\geqslant  T_3\},
$$
we derive a contradiction. This proves the theorem.
\end{proof}

\subsection{The weak limit of $\frac{u_t}{u}$}
Motivated by the exponential growth rate of $u$ obtained in the previous
section, we introduce the normalized speed 
$$
\eta(r,t) := \frac{u_t (r,t)}{u(r,t)},\quad r\in [0,1],\ t\geqslant 0.
$$
In this part we show the weak convergence of $\eta$ to $N-1$.

\begin{thm}\label{thm:weak-conv}
As $t\to\infty$, $\displaystyle \frac{u_t}{u} \stackrel{\mathrm{w}}{\to} N-1$ in $C([0,1])$.
\end{thm}

\begin{proof} 
A direct calculation shows that 
\begin{equation}\label{P-eta}
 \left\{
 \begin{array}{ll}
 \eta_t = \displaystyle \frac{\eta_{rr}}{1+u^2_r} + b_1 (r,t)\eta_r + c_1(r,t,\eta)\eta,& r\in (0,1),\ t>0,\\
 \eta_r(0,t)=0,\quad \eta_r(1,t)=0, & t>0,\\
 \eta(r,0)=\displaystyle  \eta_0(r):= \frac{u''_0}{u_0 [1+(u'_0)^2]} + \frac{(N-1)u'_0}{r u_0}, & r\in (0,1].
 \end{array}
 \right.
\end{equation}
where
$$
 b_1(r,t):=  \frac{2u_r}{u(1+u^2_r)} -
 \frac{2u_r u_{rr}}{(1+u_r^2)^2 } + \frac{N-1}{r},
$$
\begin{eqnarray*}
c_1(r,t,\eta) & := & \frac{u_{rr}}{u(1+u_r^2)} -\frac{2u^2_r u_{rr}}{u(1+u_r^2)^2} +\frac{(N-1)u_r}{ru} -\eta\\
& = & \frac{2u_r^2 }{1+u_r^2} \left[ \frac{(N-1)u_r}{ru} - \eta\right].
\end{eqnarray*}
By Corollary \ref{cor:u_r/ru}, $\frac{u_r}{ru}$ is positive and has upper bound $K_0$ in $(0,1]$. If we take $\eta= \overline{\eta}:= (N-1)K_0 + \|\eta_0\|_{L^\infty}$, then $c_1 (r,t,\bar{\eta})\leqslant 0$.
This implies that $\overline{\eta}$ is a supersolution of \eqref{P-eta}, and so $\eta$ is bounded:
\begin{equation}\label{lower-upper-eta}
0< \eta(r,t) \leqslant \overline{\eta}, \quad r\in [0,1],\ t\gg 1.
\end{equation}
Here, the first inequality follows from Theorem \ref{ut_conv}.

On the other hand,
$$
\eta = \frac{u_{rr}}{(1+u_r^2)u} + \frac{(N-1)u_r}{ru}
= \frac{(\arctan u_r)_r}{u}+\frac{(N-1)u_r}{ru}.
$$
For any test function $\rho(r) \in C^1([0,1])$ with $\rho$ and $\rho_r$ being
bounded in $[0,1]$, multiplying the above equality by $\rho$ and integrating over $[0,1]$ we have
\begin{equation}\label{Hest}
\int_0^1\Big[\eta - \frac{(N-1)u_r}{ru}\Big]\rho\,dr = (\arctan u_r)\frac{\rho}{u}\Big|_0^1 - \int_0^1\arctan u_r\frac{\rho_ru-\rho u_{r}}{u^2}\,dr.
\end{equation}
Since $u\to \infty\ (t\to \infty)$ we see that the righthand side tend to $0$, so does the lefthand
side. Recalling
$$
0\leqslant \frac{u_r}{ru} =\frac{\zeta}{r}\leqslant K_0,\qquad r\in (0,1],\ t\gg 1,
$$
and using the estimate in Theorem \ref{mainthm1} (i), we see that
\begin{equation}\label{weak-conv}
\lim\limits_{t\to \infty} \int_0^1 [\eta (r,t)-(N-1)]\rho(r) =0.
\end{equation}
This limit holds for any test function $\rho\in C^1([0,1])$, so does for any $\rho$ being taken from the space of bounded variation functions $BV([0,1])$. Thus, $\eta(\cdot,t)$ converges as $t\to \infty$ weakly to $N-1$ in $C([0,1])$.
This completes the proof of the theorem.
\end{proof}

Our results show that, the instantaneous speed has the same exponential growth rate as
the solution itself. This confirms that the higher-dimensional Robin problem exhibits an
{\it accelerating propagation mechanism}, in sharp contrast with the finite-speed
translation observed in the planar case.

\section{Asymptotic Behaviour III: The Asymptotic Mean Curvature}

In the previous sections, we have shown that both the height of the graph and
the instantaneous speed become unbounded as $t\to\infty$.
In this final asymptotic analysis, we turn to the geometric quantity of the
flow, namely the mean curvature.

The goal of this section is to show that, despite the degenerating behaviour 
of the graph itself, the mean curvature admits a well-defined weak asymptotic
limit. More precisely, we assume $N\geqslant 2$ and  show that the mean curvature $H$ of the graph $(x,u(x,t))$ converges in certain weak sense to $\frac{N-1}{r}$. 

We will use the notation in the previous sections.
First, as a consequence of the last theorem in the previous section we have the following weak convergence for $rH$.

\begin{prop}\label{prop:rH-conv}
$r H(r,t) \stackrel{\mathrm{w}}{\to} N-1\ (t\to\infty)$ in $C([0,1])$.
\end{prop}

\begin{proof}
Note that
$$
\eta = \frac{u_t}{u} = \frac{u_{rr}}{u(1+u_r^2)} +\frac{(N-1)u_r}{ru} = H \frac{\sqrt{1 +u_r^2}}{u}. $$
Our conclusion follows directly from Theorem \ref{thm:weak-conv} if we can show that
$$
\left\| \frac{\sqrt{1 +u_r^2}}{u} - r\right\|_{C([0,1])}\to 0 \mbox{\ \ as \ \ }t\to \infty.
$$
This limit can be shown easily. In fact, for any small $\varepsilon>0$, there exists $T_0>0$ large such that when $t\geqslant T_0$,
$$
\frac{1}{u^2}\leqslant \varepsilon^2,\quad \left| \frac{u_r}{u}-r\right| \leqslant C e^{-\frac{N-1}{3}t} \leqslant \varepsilon^2.
$$
Then, on $[0,\varepsilon]$ we have
$$
\left| \frac{\sqrt{1 +u_r^2}}{u} - r\right| \leqslant \left| \sqrt{\frac{1}{u^2} +\frac{u_r^2}{u^2}} +r\right| \leqslant \left|\sqrt{\varepsilon^2 + (r+\varepsilon^2)^2} +r\right| < (\sqrt{5}+1)\varepsilon.
$$
On $[\varepsilon,1]$ we have
$$
\left| \frac{\sqrt{1 +u_r^2}}{u} - r\right|  =  \left| \frac{ \frac{1}{u^2}  + \frac{u_r^2}{u^2} -r^2} {\sqrt{\frac{1}{u^2} + \frac{u_r^2}{u^2}} +r}  \right|  \leqslant  \frac{\varepsilon^2 + \left| \frac{u_r^2}{u^2} - r^2\right|}{r}  \leqslant \frac{\varepsilon^2 +  (2r+\ve^2) \ve^2 }{\ve} <4 \varepsilon.
$$
This proves the proposition.
\end{proof}

Next we present another convergence of $H$ to $\frac{N-1}{r}$ as shown in the following result.

\begin{thm}\label{thm:L2-conv}
For any given $\ve\in(0,\frac{1}{6})$, there holds
$$
\int_0^\infty e^{2(N-1-3\ve)t}\int_{e^{-\ve t}}^1\Big|H(r,t) - \frac{N-1}{r}\Big|^2\,drdt < \infty.
$$
\end{thm}

\begin{proof}
For given $\varepsilon \in (0,\frac16)$, we first show
$$
\int_0^\infty e^{2(N-1-3\ve)t}\int_{e^{-\ve t}}^1\Big|\frac{u_r}{r\sqrt{1+u_r^2}}- \frac{1}{r}\Big|^2\,drdt <\infty.
$$
When $t$ is large, say $t> T_1$, by Lemma \ref{lem:zeta-lower-est} there holds
$$
\zeta = w_r \geqslant   r -   \frac{1}{2}e^{-\ve t}
\geqslant \frac{1}{2}e^{-\ve t},\qquad r\geqslant  e^{-\ve t},\ t>T_1.
$$
So we have $u_r =  u w_r \geqslant \frac12 e^{(N-1-\ve)t}$ for $r\geqslant e^{-\ve t},\ t> T_1$. Hence
\begin{eqnarray*}
\displaystyle \int_{e^{-\ve t}}^1 \Big|\frac{u_r}{r\sqrt{1+u_r^2}} - \frac{1}{r}\Big|^2\,dr
& = & \displaystyle \int_{e^{-\ve t}}^1 \frac{dr}{r^2(1+u_r^2)[\sqrt{1+u_r^2}+u_r]^2} \leqslant
\int_{e^{-\ve t}}^1 \frac{dr}{4r^2 u_r^4}\\
&\leqslant & \displaystyle    \frac{4}{ e^{-2\ve t}e^{4(N-1-\ve)t}} = 4 e^{-4(N-1-\frac{3}{2}\ve)t},\qquad t>T_1.
\end{eqnarray*}
Consequently,
\begin{equation}\label{meanc_1}
\int_0^\infty e^{4(N-1-2\ve)t}\int_{e^{-\ve t}}^1\Big|\frac{u_r}{r\sqrt{1+u_r^2}}- \frac{1}{r}\Big|^2\,drdt<\infty.
\end{equation}

Next, we estimate the integral
$$
\int_0^\infty e^{2(N-1-3\ve)t}\int_{e^{-\ve t}}^1\Big|\frac{u_{rr}}{(1+u_r^2)^{3/2}}\Big|^2\,drdt.
$$
Using the equation of $u$ we have
\begin{eqnarray*}
0 & = & \displaystyle \int_{e^{-\varepsilon t}}^1 2u_{rr}\Big[-u_t + \frac{u_{rr}}{1+u_r^2}+\frac{N-1}{r}u_r\Big]\,dr\\
& =  & \displaystyle  \Big(-2u_ru_t+\frac{N-1}{r}u_r^2\Big)\Big|_{e^{-\varepsilon t}}^1 +  \int_{e^{-\varepsilon t}}^1\Big[2u_ru_{rt}+\frac{2u_{rr}^2}{1+u_r^2}+\frac{N-1}{r^2}u_r^2\Big]\,dr\\
& =  &   I_1 + \frac{d}{dt} \left[ \int_{e^{-\varepsilon t}}^1 u_r^2 dr\right]  +  \int_{e^{-\varepsilon t}}^1 \frac{2u_{rr}^2}{1+u_r^2} \,dr,
\end{eqnarray*}
with
\begin{eqnarray*}
I_1 & := & -2u(1,t) u_t(1,t) + (N-1)u^2(1,t) +2u_r({e^{-\varepsilon t}},t)u_t ({e^{-\varepsilon t}},t) - (N-1){e^{\varepsilon t}} u_r^2 ({e^{-\varepsilon t}},t)\\
& & - \varepsilon {e^{-\varepsilon t}} u_r^2 ({e^{-\varepsilon t}},t) + (N-1)\int_{{e^{-\varepsilon t}}}^1 \frac{u_r^2}{r^2} dr.
\end{eqnarray*}
Set $\beta := 2(N-1+\varepsilon)$ and $I_2 := I_1 + \beta \int_{e^{-\varepsilon t}}^1 u_r^2 dr$. Then \begin{equation}\label{3-terms}
0=  e^{-\beta t} I_2 + \frac{d}{dt} \left[ e^{-\beta t} \int_{e^{-\varepsilon t}}^1 u_r^2 dr \right] +  e^{-\beta t} \int_{e^{-\varepsilon t}}^1 \frac{2u_{rr}^2}{1+u_r^2} \,dr.
\end{equation}
By Corollary \ref{cor:C0est-v-u}, Lemma \ref{lem:zeta-upper-est} and \eqref{lower-upper-eta} we see that
$$
u(r,t),\ \ |u_r(r,t)|,\ \ |u_t(r,t)|\leqslant C e^{(N-1) t},\quad r\in [0,1],\ t>0,
$$
for some $C>0$. Hence the absolute value of each term in $I_2$ is less than $C_1 e^{2(N-1)t + \varepsilon t}$ for some $C_1>0$. It follows that, for some $C_2>0$,
$$
\int_0^\infty \Big| e^{-\beta t} I_2 (t) \Big| dt \leqslant C_2 \int_0^\infty e^{[-\beta  + 2(N-1)+\varepsilon] t}  dt <\infty.
$$
On the other hand,
\begin{eqnarray*}
\int_0^\infty \frac{d}{dt} \left[ e^{-\beta t} \int_{e^{-\varepsilon t}}^1 u_r^2 dr \right] dt & =&
\lim\limits_{T\to \infty} \int_0^T \frac{d}{dt} \left[ e^{-\beta t} \int_{e^{-\varepsilon t}}^1 u_r^2 dr \right] dt \\
& = & \lim\limits_{T\to \infty} e^{-\beta T} \int_{e^{-\ve T}}^1 u_r^2 (r,T)dr =0.
\end{eqnarray*}
Thus, by \eqref{3-terms} we have
\begin{equation}\label{e-beta-t}
\int_0^\infty e^{-\beta t} \int_{e^{-\varepsilon t}}^1\frac{u_{rr}^2}{1+u_r^2}\,drdt < \infty.
\end{equation}
Finally, using $u_r\geqslant  \frac12 e^{(N-1-\ve)t}$ for $r \geqslant e^{-\ve t}$ we obtain
$$
\int_{e^{-\ve t}}^1\left[ \frac{u_{rr}}{(1+u_r^2)^{3/2}} \right]^2\,dr \leqslant
\int_{e^{-\ve t}}^1\frac{u_{rr}^2}{1+u_r^2}\cdot\frac{dr}{u_r^4} \leqslant
16 e^{-4(N-1-\ve)t}\int_{e^{-\ve t}}^1\frac{u_{rr}^2}{1+u_r^2}\,dr.
$$
It follows from \eqref{e-beta-t} that
\begin{equation}\label{meanc_2}
\int_0^\infty e^{2(N-1-3\ve)t}\int_{e^{-\ve t}}^1\frac{u_{rr}^2}{(1+u_r^2)^{3}}\,drdt \leqslant  16 \int_0^\infty e^{-\beta t}\int_{e^{-\ve t}}^1\frac{u_{rr}^2}{1+u_r^2}drdt  < \infty.
\end{equation}
Combining \eqref{meanc_1} with \eqref{meanc_2} we obtain the desired result.
\end{proof}

\begin{remark}\rm
This completes the description of the long-time dynamics.
Although the solution profile degenerates and the propagation speed diverges,
the mean curvature converges to a stable geometric limit.
\end{remark}

\noindent
{\textbf{Proof of Main Theorem}}. The conclusion (i) follows from Lemma \ref{lem:lower-est-u}, Corollary \ref{cor:C0est-v-u}, Theorem \ref{ut_conv} and Lemma \ref{lem:lower-est-ur}. The conclusions in (ii)
follow from Theorem \ref{thm:weak-conv} and Proposition \ref{prop:rH-conv}. These conclusions hold for $N\geqslant 2$. When $N\geqslant 3$, the conclusions (iii) and (iv) are proved in Theorem \ref{mainthm1}.
\qed

\medskip
\noindent
{\bf Acknowledgements}. The second author would like to thank Professors Xinan Ma, Xujia Wang and Xingping Zhu for helpful discussions.

\medskip
\noindent
{\bf Conflict of Interest}. 
The authors declare that they have no conflict of interest.

\medskip
\noindent
{\bf Data Availability}. 
All results in this paper are derived from mathematical proofs and do not rely on external datasets.

%
%
%
%
%
%
%
%

\section*{Appendix}
In this part we give a formal derivation of the mean curvature flow equation \eqref{mcf1} and the Robin boundary condition \eqref{b3} from singular limit problem of the Allen-Cahn equation.

Let $\widetilde{\Omega}\subset \R^{N+1}$ be a domain with smooth boundary, $\tilde{\nu}$ be the unit inner normal vector of $\partial \widetilde{\Omega}$, $\varepsilon>0$ be a small parameter and $\tilde{u}^\varepsilon (\tilde{x},t)$  be a classical solution of the following problem:
\begin{equation}\label{AC-p}
\left\{
\begin{array}{ll}
\displaystyle \tilde{u}_t = \Delta \tilde{u} + \frac{1}{\varepsilon^2} [\tilde{u}-\tilde{u}^3], & \tilde{x}\in \widetilde{\Omega}, t>0,\\
\displaystyle \frac{\partial \tilde{u}}{\partial \tilde{\nu}} =\frac{1}{\varepsilon} \tilde{g} (\tilde{x},t,\tilde{u}), & \tilde{x}\in \partial \widetilde{\Omega}, t>0,\\
\tilde{u}(x,0)=\xi(\tilde{x}), & \tilde{x}\in \overline{\widetilde{\Omega}},
\end{array}
\right.
\end{equation}
where $\tilde{g}$ and $\xi$ are smooth functions with $\tilde{g}(x,t,\pm 1)\equiv 0$. For any $t\geqslant 0$, denote $\Gamma^\varepsilon(t):= \{x\in \widetilde{\Omega} \mid \tilde{u}^\varepsilon (\tilde{x},t)=0\}$ and denote $d^\varepsilon(\tilde{x},t)$ the signed distance function from $\tilde{x}\in \widetilde{\Omega}$ to $\Gamma^\varepsilon (t)$.

In the case where $\Gamma^\varepsilon(t)$ is a family of simple closed hypersurfaces in $\widetilde{\Omega}$, it was shown in  \cite{Matano2008, Chen1992, Fei, NMHS} etc., rigorously or formally, that as $\varepsilon\to 0$,  $\Gamma^\varepsilon(t)\to \Gamma(t)$, $d^\varepsilon(\tilde{x},t)\to d(\tilde{x},t)$, and $d$ satisfies
$d_t =\Delta_{\tilde{x}} d$ on $\Gamma(t)$. This equation is nothing but the mean curvature flow since $d_t$ and $\Delta_{\tilde{x}} d$ equal to the normal velocity and the mean curvature of $\Gamma(t)$, respectively. In particular, if $\Gamma(t)$ is the graph of a function $x_{N+1}=u(x_1,\cdots, x_N,t)$, then the equation $d_t=\Delta_{\tilde{x}} d$ is converted into \eqref{mcf1}.

We now consider the case where $\Gamma^\varepsilon(t)$ contacts $\partial \widetilde{\Omega}$. We will use the matched asymptotic expansion method as in \cite{Fei, NMHS} to derive formally the boundary condition satisfied by $\Gamma(t)$. More precisely, assume, for each $t>0$, $\Gamma^\varepsilon (t)$ is a simple $N$-dimensional hypersurface contacting $\partial \widetilde{\Omega}$ at a simple $(N-1)$-dimensional hypersurface $\gamma(t)$. So $\Gamma^\varepsilon(t)$ separate $\widetilde{\Omega}$ into two connect components, denote them by $\widetilde{\Omega}_\pm (t)$ respectively.
For $\tilde{x}$ lying far from $\Gamma^\varepsilon(t)$ we substitute the {\it outer expansion}: $\tilde{u}^\varepsilon(\tilde{x},t) =  \tilde{u}_0(\tilde{x},t) + \varepsilon \tilde{u}_1(\tilde{x},t) + \varepsilon^2 \tilde{u}_2 (\tilde{x},t)+\cdots$ into \eqref{AC-p} and collect the terms of $\varepsilon^{-2},\varepsilon^{-1},\cdots$, to conclude that
$$
\tilde{u}_0(\tilde{x},t) = \pm 1 \mbox{ for }\tilde{x}\in \widetilde{\Omega}_\pm (t) \mbox{ and }|d^\varepsilon(\tilde{x},t)|\gg \varepsilon,\qquad  \tilde{u}_i(\tilde{x},t) \equiv 0\ (i=1,2,\cdots).
$$
Next, we define $y:= \frac{d^\varepsilon(\tilde{x},t)}{\varepsilon} =\frac{1}{\varepsilon} \left[ d_0(\tilde{x},t) + \varepsilon d_1(\tilde{x},t) +\cdots\right]$
and suppose that $\tilde{u}^\varepsilon$ has the following {\it inner expansion} near $\Gamma^\varepsilon (t)$:
$$
\tilde{u}^\varepsilon (\tilde{x},t) = U_0(\tilde{x},y,t) + \varepsilon U_1(\tilde{x},y,t) + \varepsilon^2
U_2(\tilde{x},y,t) + \cdots.
$$
Substituting it into \eqref{AC-p}, collecting the terms of $\varepsilon^{-2}, \varepsilon^{-1}$ and
using the matching condition and normalized conditions (see for example \cite{NMHS}) we have
\begin{equation}\label{U0-p}
\left\{
\begin{array}{l}
U_{0yy}+ U_0 -U_0^3 =0, \qquad \tilde{x}\in \Gamma(t),\ y\in \R,\ t>0,\\
U_0(\tilde{x},-\infty,t)=1,\qquad U_0(\tilde{x},\infty,t)=-1,\qquad U_0(\tilde{x},0,t)=0, \qquad  \tilde{x}\in \Gamma(t),\ t>0,
\end{array}
\right.
\end{equation}
\begin{equation}\label{U1-p}
\left\{
\begin{array}{ll}
U_{1yy}+ [1 -3 U_0^2]U_1 = (d_{0t}-\Delta_{\tilde{x}}d_0) U_{0y}, & \tilde{x}\in \Gamma(t),\ y\in \R,\ t>0,\\
U_1(\tilde{x},-\infty,t)=U_1(\tilde{x},\infty,t)=U_1(\tilde{x},0,t)=0, & \tilde{x}\in \Gamma(t),\ t>0,
\end{array}
\right.
\end{equation}
and
\begin{equation}\label{d0-bdry}
U_{0y} D d_0 \cdot \tilde{\nu}=\tilde{g}(\tilde{x},t,U_0(\tilde{x},0,t)), \qquad \tilde{x}\in \gamma(t), \  t>0.
\end{equation}
By \eqref{U0-p} we can choose $U_0 (\tilde{x},y,t)\equiv \varphi(y):= \tanh \big( \frac{-y}{\sqrt{2}}\big)$. Differentiating the equation in \eqref{U0-p} with respect to $y$ we have $[\varphi']''+[1-3U_0^2]\varphi'=0$. This means that $\varphi'(y)$ is the eigenfunction of the operator $\partial^2_{y} + [1-3U_0^2 ]I$ corresponding to eigenvalue $0$. By the Fredholm theorem, \eqref{U1-p} has a solution iff
$$
0= \int_\R (d_{0t}-\Delta_{\tilde{x}}d_0) [\varphi'(y)]^2 dy = (d_{0t}-\Delta_{\tilde{x}}d_0) \cdot \int_\R [\varphi'(y)]^2 dy.
$$
Hence we obtain $d_{0t}=\Delta_{\tilde{x}} d_0$ for $\tilde{x}\in \Gamma(t)$, which, as we have mentioned above, is a mean curvature flow equation, and it is converted into \eqref{mcf1} when $\Gamma(t)$ is the graph of a function $x_{N+1}=u(x_1,\cdots, x_N,t)$.
Finally, by \eqref{d0-bdry} we have
\begin{equation}\label{bdry-00}
 D d_0(\tilde{x},t)\cdot \tilde{\nu} = \frac{\tilde{g}(\tilde{x},t,0)}{\varphi'(0)}= -\sqrt{2}\tilde{g}(\tilde{x},t,0),
\qquad \tilde{x}\in \gamma(t),\ t>0.
\end{equation}
In the special case where $\widetilde{\Omega}$ is a cylinder $\Omega\times \R$ for some bounded $\Omega\subset \R^N$ and $\Gamma(t)$ is the graph of $x_{N+1}=u(x,t)$ for $x=(x_1,\cdots, x_N)\in \Omega$, we see that $\gamma(t)$ is a $(N-1)$-dimensional hypersurface on $\partial\Omega \times \R$. In addition, $\tilde{\nu}=(\nu,0)$, where $\nu$ is the inner unit normal vector of $\partial \Omega$. Assume the signed distance function $d^\varepsilon(\tilde{x},t)$ is positive when $\tilde{x}$ lies above $\Gamma(t)$, then
$$
D d_0(x,x_{N+1},t) = \frac{(-D u, 1)}{\sqrt{1+| D u|^2}},
$$
and so \eqref{bdry-00} can be rewritten as
\begin{equation}\label{MCF-BC}
\frac{ -D u\cdot \nu }{\sqrt{1+|D u|^2}} = g(x,t,u(x,t)) := -\sqrt{2}\tilde{g}(x,u(x,t),t,0),\qquad x\in \partial \Omega,\ t>0.
\end{equation}
This leads to \eqref{b3}.

\end{document}